\documentclass[a4paper,12pt]{elsarticle}
\usepackage[utf8]{inputenc}
\usepackage{amssymb}
\usepackage{amsfonts}
\usepackage{amsmath,amsthm}

\usepackage{enumerate}

\newtheorem{theorem}{Theorem}
\newtheorem{observation}[theorem]{Observation}
\newtheorem{conjecture}[theorem]{Conjecture}

\newtheorem{corollary}[theorem]{Corollary}
\newtheorem{lemma}[theorem]{Lemma}

\newcommand{\Aut}{\ensuremath{\operatorname{Aut}}}
\newcommand{\id}{\ensuremath{\operatorname{id}}}
\newcommand{\orG}{\ensuremath{\overrightarrow{G}}}

\begin{document}

\title{Arc-distinguishing of orientations of graphs}

\author[AGH]{Aleksandra Gorzkowska} 
\ead{agorzkow@agh.edu.pl}

\author[AGH]{Jakub Kwa\'sny\corref{cor1}} 
\ead{jkwasny@agh.edu.pl}

\address[AGH]{AGH University of Krakow, al. A. Mickiewicza 30, 30-059 Krakow, Poland}
\cortext[cor1]{Corresponding author}

\begin{abstract}
The distinguishing index $D'(G)$ of a graph $G$ is the minimum number of colours in an edge colouring preserved only by the identity automorphism. We study how orienting the edges affects this parameter, relating the minimum and maximum distinguishing indices over all orientations of $G$ to $D'(G)$. We establish bounds and exact relations for bipartite graphs and trees, and study rigid orientations of traceable and claw-free graphs. Our results answer a question of Meslem and Sopena on unbalanced complete bipartite graphs.

\end{abstract}


\maketitle

\section{Introduction}

Edge colours and edge directions provide two ways to break graph symmetries. We study how they interact: how much can orienting the edges reduce the number of colours needed to distinguish a graph?

We follow the terminology and notation of \cite{west}, and write $[k]=\{1,\ldots,k\}$. All graphs considered are connected and distinct from $K_2$. Edge colourings need not be proper. A colouring $c\colon E(G)\to[k]$ \emph{breaks} an automorphism $\varphi\in\Aut(G)$ if some edge $xy$ satisfies $c(\varphi(x)\varphi(y))\neq c(xy)$. The colouring is \emph{distinguishing} if it breaks every non-trivial automorphism. The \emph{distinguishing index} $D'(G)$ is the minimum number of colours in such a colouring. We exclude $K_2$ because the automorphism interchanging its vertices preserves its only edge, regardless of its colour.

Distinguishing vertex colourings were studied by Babai \cite{Babai} in 1977 and later by Albertson and Collins \cite{AlbertsonCollins} in 1996. Kalinowski and Pil\'sniak \cite{KalinowskiPilsniak} introduced the distinguishing index in 2015. Subsequent work includes bounds for traceable graphs \cite{Pilsniak_traceable}, claw-free graphs \cite{GorzkowskaKargul} and regular graphs \cite{KwasnyStawiski_regular}, as well as bounds in terms of the maximum degree and a classification of the equality cases \cite{Pilsniak_traceable}.

An automorphism of a digraph preserves arcs and their directions. A distinguishing arc colouring and the distinguishing index of a digraph are therefore defined in the same way, with arcs in place of edges. Related work concerns symmetric digraphs, obtained by replacing each edge of a graph with two opposite arcs \cite{KalinowskiPilsniak2,KalinowskiPilsniakProrok}. Here we consider \emph{orientations}: digraphs $\orG$ obtained from $G$ by replacing each edge $xy$ with exactly one of the arcs $\overrightarrow{xy}$ and $\overrightarrow{yx}$.

Meslem and Sopena \cite{MeslemSopena} introduced the parameters
\[
OD'^-(G)=\min D'(\orG)
\qquad\text{and}\qquad
OD'^+(G)=\max D'(\orG),
\]
where both extrema range over all orientations of $G$. They determined these parameters for paths, cycles, complete graphs and balanced complete bipartite graphs. We extend the study to broader classes, focusing on the relation between $OD'^-(G)$, $OD'^+(G)$ and $D'(G)$ rather than on individual values. An orientation is called \emph{rigid} if its only automorphism is the identity; thus $OD'^-(G)=1$ precisely when $G$ admits a rigid orientation.

Orienting edges can only remove automorphisms, thus the following observation which is used throughout the paper.

\begin{observation} \label{obs: equal_aut}
Let $\orG$ be an orientation of a graph $G$. Then:
\begin{enumerate}[(i)]
\item $\Aut(\orG) \subseteq \Aut(G)$,
\item if $\Aut(\orG) = \Aut(G)$, then $D'(\orG) = D'(G)$,
\item if $\Aut(\orG) = \{\id\}$, then $D'(\orG) = 1$.
\end{enumerate}
\end{observation}

For a graph or digraph $H$, a vertex set $S$ is \emph{setwise fixed} if $\varphi(S)=S$ for every $\varphi\in\Aut(H)$, and \emph{pointwise fixed} if $\varphi(v)=v$ for every $v\in S$ and every $\varphi\in\Aut(H)$. In particular, a setwise fixed set may have its vertices permuted, whereas a pointwise fixed set may not. A vertex is \emph{fixed} if every automorphism fixes it.

In Section \ref{sec:bipartite}, we determine $OD'^-(G)$ and $OD'^+(G)$ in terms of $D'(G)$ for bipartite graphs whose bipartition classes cannot be interchanged by an automorphism. This answers a question of Meslem and Sopena. For arbitrary connected bipartite graphs, we restrict each parameter to at most two consecutive integer values. We then give exact formulas for trees, with the cases distinguished by the number of non-isomorphic optimal colourings of the relevant rooted subtrees.

In Section \ref{sec:D'=2}, we construct rigid orientations of traceable graphs and of connected claw-free graphs of order at least six. We also give a criterion for $OD'^+(G)=1$ when $D'(G)=2$. As a tool for the claw-free case, we prove that every Hamiltonian graph other than a cycle admits a strongly connected rigid orientation in which a prescribed Hamiltonian cycle is the unique directed Hamiltonian cycle.

\section{Bipartite graphs} \label{sec:bipartite}

We begin with the part of the result of Meslem and Sopena \cite{MeslemSopena} concerning distinguishing indices of orientations of unbalanced complete bipartite graphs.

\begin{theorem} \cite{MeslemSopena}
For every two integers $m$ and $n$, $2 \leq m < n$, the following hold:
\begin{enumerate}
    \item $OD'^+(K_{m,n})=D'(K_{m,n})$.
    \item If $K_{m,n}$ admits a rigid orientation, then $OD'^-(K_{m,n})=1$.
    \item If $K_{m,n}$ does not admit any rigid orientation, then $OD'^-(K_{m,n}) \leq D'(K_{m,\lceil \frac{n}{m-1}\rceil})$.
\end{enumerate}
\end{theorem}

To extend these results beyond complete bipartite graphs, we use the preservation of the bipartition classes by automorphisms. The following theorem isolates this property for partitions into any number of independent sets.

\begin{theorem} \label{lemma: part_indep}
Let $G=(V,E)$ be a graph with a partition $V=V_1\cup\dots\cup V_k$ into $k\geq 1$ independent sets, each setwise fixed by every automorphism of $G$. Then $OD'^+(G)=D'(G)$ and $OD'^-(G)=\lceil D'(G)/2\rceil$.
\end{theorem}

\begin{proof}
For $OD'^+(G)$, orient every edge from $V_i$ to $V_j$ whenever $i<j$. Since each $V_i$ is independent and setwise fixed, every automorphism of $G$ preserves this orientation. Thus $\Aut(\orG)=\Aut(G)$, and Observation \ref{obs: equal_aut} gives $D'(\orG)=D'(G)$. Since a distinguishing edge colouring of $G$ also distinguishes every orientation of $G$, this proves $OD'^+(G)=D'(G)$.

For $OD'^-(G)$, we encode edge colours by arc colours and directions. Let $r\geq 1$ and consider an edge colouring $c=(c_1,c_2)\colon E\to\{0,1\}\times[r]$. For each edge $uv$ with $u\in V_i$, $v\in V_j$ and $i<j$, orient it from $u$ to $v$ if $c_1(uv)=0$, and from $v$ to $u$ otherwise. Give the resulting arc colour $c_2(uv)$. Conversely, every orientation with an arc colouring from $[r]$ determines such an edge colouring.

Because every automorphism of $G$ preserves each $V_i$, it preserves $c_1$ if and only if it preserves the orientation. It therefore preserves $c$ if and only if it is an automorphism of $\orG$ preserving $c_2$. Hence $c$ is distinguishing on $G$ if and only if $c_2$ is distinguishing on $\orG$.

Taking $r=\lceil D'(G)/2\rceil$, we can encode a distinguishing edge colouring of $G$ in this way and obtain an orientation with distinguishing index at most $r$. Thus $OD'^-(G)\leq\lceil D'(G)/2\rceil$. Conversely, an orientation with a distinguishing arc colouring using $r$ colours yields a distinguishing edge colouring of $G$ with at most $2r$ colours. Therefore $D'(G)\leq 2r$, which gives the reverse inequality.
\end{proof}

For bipartite graphs, it suffices to exclude automorphisms that interchange the two classes.

\begin{corollary} \label{cor: bipartite}
Let $G=(X\cup Y,E)$ be a bipartite graph with no automorphism interchanging $X$ and $Y$. Then $OD'^-(G)=\lceil D'(G)/2\rceil$ and $OD'^+(G)=D'(G)$.
\end{corollary}
\begin{proof}
Apply Theorem \ref{lemma: part_indep} with $V_1=X$ and $V_2=Y$.
\end{proof}

For $m<n$, the two classes of $K_{m,n}$ have different sizes and cannot be interchanged. Corollary \ref{cor: bipartite} therefore answers the question of Meslem and Sopena \cite{MeslemSopena} about $OD'^-(K_{m,n})$ when $n$ is substantially larger than $m$. To express the answer in terms of the class sizes, we use the following result of Fisher and Isaak \cite{FisherIsaak} and Imrich, Jerebic and Klav{\v{z}}ar \cite{ijk}.

\begin{theorem} \cite{FisherIsaak,ijk}
Let $m,n,r$ be integers such that $2\leq m<n$, $r\geq 2$ and $(r-1)^m<n\leq r^m$. Then
$$D'(K_{m,n}) = \left\{ 
\begin{array}{ll}
r, & \textrm{if } n \leq r^m - \lceil \log_r m \rceil -1;\\
r+1, & \textrm{if } n \geq  r^m - \lceil \log_r m \rceil +1.
\end{array} 
\right.$$
If $n=r^m - \lceil \log_r m \rceil$, then $D'(K_{m,n})$ is either $r$ or $r+1$ and can be computed recursively in time $O(\log^*(n))$.
\end{theorem}

Combining this result with Corollary \ref{cor: bipartite} gives the following formula for $2\leq m<n$. To include the balanced case, we anticipate Theorem \ref{thm: traceable}, proved in Section \ref{sec:D'=2}: every $K_{m,m}$ has a Hamiltonian path, so $OD'^-(K_{m,m})=1$, in agreement with the formula below.

\begin{corollary}
Let $m,n,r$ be integers such that $2\leq m\leq n$, $r\geq 2$ and $(r-1)^m<n\leq r^m$. Then
$$OD'^-(K_{m,n}) = \left\{ 
\begin{array}{ll}
\lceil \frac{r}2 \rceil, & \textrm{if } n \leq r^m - \lceil \log_r m \rceil -1;\\
\lceil \frac{r+1}2 \rceil, & \textrm{if } n \geq  r^m - \lceil \log_r m \rceil +1.
\end{array} 
\right.$$
If $n=r^m - \lceil \log_r m \rceil$, then $OD'^-(K_{m,n})$ is either $\lceil r/2 \rceil$ or $\lceil (r+1)/2 \rceil$ and can be computed recursively in time $O(\log^*(n))$.
\end{corollary}

In particular, an unbalanced complete bipartite graph $K_{m,n}$ admits a rigid orientation if and only if $D'(K_{m,n}) = 2$. 

We now allow automorphisms that interchange the bipartition classes. Although the equalities in Corollary \ref{cor: bipartite} need not hold, each parameter is still restricted to at most two consecutive integer values in terms of $D'(G)$.

\begin{theorem} \label{thm: bipartite_bounds}
Let $G$ be a connected bipartite graph of order at least three. Then
\[
\left\lfloor\frac{D'(G)}{2}\right\rfloor
\leq OD'^-(G)
\leq \left\lfloor\frac{D'(G)}{2}\right\rfloor+1
\]
and
\[
D'(G)-1\leq OD'^+(G)\leq D'(G).
\]
\end{theorem}
\begin{proof}
Let $G=(X\cup Y,E)$. Since $G$ is connected, every automorphism either preserves both bipartition classes setwise or interchanges them.

We first note that any edge colouring with $t$ colours that breaks all non-trivial automorphisms preserving both classes implies $D'(G)\leq t+1$. Indeed, at most one non-trivial automorphism can preserve such a colouring. Otherwise, two distinct colour-preserving automorphisms $\varphi$ and $\psi$ would both interchange the classes, so $\varphi^{-1}\psi$ would be a non-trivial colour-preserving automorphism preserving both classes, a contradiction. If a non-trivial colour-preserving automorphism remains, it moves some edge $e$, since $G$ is connected and has at least three vertices. Recolour $e$ with a new colour. Every automorphism preserving the new colouring fixes $e$ setwise and also preserves the original colouring, so it must be the identity. If the original colouring is already distinguishing, no additional colour is needed.

To prove the lower bound for $OD'^-(G)$, take any orientation $\orG$ with a distinguishing arc colouring using $k$ colours. Apply the same encoding of directions and arc colours as in the proof of Theorem \ref{lemma: part_indep}, taking $V_1=X$ and $V_2=Y$. This gives an edge colouring with at most $2k$ colours that breaks all non-trivial automorphisms preserving both bipartition classes: any such automorphism preserving the encoded colouring would also preserve the orientation and its distinguishing arc colouring. The preceding observation gives $D'(G)\leq 2k+1$. Taking $k=OD'^-(G)$ yields the desired lower bound.

For the upper bound, put $d=D'(G)$ and take a distinguishing edge colouring with colours $1,\ldots,d$. Keep colour $d$ in a singleton group, and partition the remaining colours into pairs, with one additional singleton if necessary. Replace the colours in each group by a common arc colour. For each pair, orient the edges of its first colour from $X$ to $Y$ and those of its second colour from $Y$ to $X$. Orient the edges in every singleton group from $X$ to $Y$. This uses
\[
1+\left\lceil\frac{d-1}{2}\right\rceil
=\left\lfloor\frac{d}{2}\right\rfloor+1
\]
arc colours. An automorphism preserving this arc colouring cannot interchange $X$ and $Y$, since the non-empty class corresponding to the singleton group $\{d\}$ contains only arcs from $X$ to $Y$. Once the bipartition classes are preserved, the original edge colour can be recovered from the arc colour and its direction. Thus such an automorphism preserves the original distinguishing colouring and is the identity.

Finally, orient every edge from $X$ to $Y$, obtaining an orientation $\orG_0$ whose automorphisms are precisely the automorphisms of $G$ preserving both bipartition classes. A distinguishing arc colouring of $\orG_0$ therefore breaks all non-trivial automorphisms of $G$ preserving both classes. Applying the initial observation with $t=D'(\orG_0)$ gives
\[
D'(G)-1\leq D'(\orG_0)\leq OD'^+(G).
\]
The reverse bound $OD'^+(G)\leq D'(G)$ follows because every distinguishing edge colouring of $G$ also distinguishes each of its orientations.
\end{proof}

We now specialise to trees. The \emph{centre} of a tree consists of a single vertex or the two endpoints of an edge, and is setwise fixed by every automorphism. If the centre is pointwise fixed, no automorphism can interchange the bipartition classes. We therefore obtain the following consequence of Corollary \ref{cor: bipartite}.

\begin{corollary} \label{cor: most_trees}
Let $T$ be a tree whose centre consists of a single vertex, or of an edge whose endpoints are fixed by every automorphism. Then $OD'^-(T)=\lceil D'(T)/2\rceil$ and $OD'^+(T)=D'(T)$.
\end{corollary}

It remains to consider trees whose central edge $e$ can be reversed by an automorphism. The two components of $T-e$ are then isomorphic as rooted trees, with the endpoints of $e$ as their roots. Let $(T',r)$ denote either component with this choice of root.

An automorphism of a rooted tree $(T',r)$ is an automorphism of $T'$ that fixes $r$. Its distinguishing index $D'((T',r))$ is the minimum number of colours in an edge colouring that breaks all non-trivial root-preserving automorphisms. A distinguishing colouring using this number of colours is called \emph{optimal}.

We compare colourings using the same fixed palette. Two edge colourings $c_1,c_2$ of $(T',r)$ are \emph{isomorphic} if there is a root-preserving automorphism $\varphi$ such that $c_2(xy)=c_1(\varphi(x)\varphi(y))$ for every $xy\in E(T')$; otherwise, they are \emph{non-isomorphic}. The next theorem shows that the existence of two non-isomorphic optimal colourings determines the relation between the orientation parameters of $T$ and $D'(T)$.

\begin{theorem} \label{thm: other_trees}
Let $T$ be a tree of order at least three whose central edge can be reversed by an automorphism, and let $(T',r)$ be either rooted component obtained by deleting that edge.
\begin{enumerate}[(i)]
\item If $(T',r)$ admits at least two non-isomorphic optimal distinguishing colourings, then
\[
OD'^+(T)=D'(T),\qquad
OD'^-(T)=\left\lceil\frac{D'(T)}{2}\right\rceil.
\]
\item If $(T',r)$ has a unique optimal distinguishing colouring up to a root-preserving automorphism, then
\[
OD'^+(T)=D'(T)-1,\qquad
OD'^-(T)=\left\lceil\frac{D'(T)-1}{2}\right\rceil.
\]
\end{enumerate}
\end{theorem}

\begin{proof}
Let $e$ be the central edge of $T$, and put $q=D'((T',r))$. In every orientation of $T$, the direction of $e$ fixes its endpoints individually and prevents the two components of $T-e$ from being interchanged. The automorphisms of either rooted oriented component extend to the whole orientation by fixing the other component. Thus an arc colouring distinguishes the orientation of $T$ if and only if it distinguishes both rooted oriented components.

The proof of Theorem \ref{lemma: part_indep} applies to root-preserving automorphisms of $T'$, since they preserve both bipartition classes. Consequently, the maximum and minimum distinguishing indices over orientations of $(T',r)$ are $q$ and $\lceil q/2\rceil$, respectively. Applying the corresponding constructions to both components and directing $e$ arbitrarily gives
\[
OD'^+(T)=q,\qquad OD'^-(T)=\lceil q/2\rceil.
\]

It remains to relate $q$ to $D'(T)$. Every distinguishing edge colouring of $T$ restricts to a distinguishing colouring of each rooted component, so $D'(T)\geq q$. If $(T',r)$ has two non-isomorphic optimal colourings, use one on each component and colour $e$ arbitrarily with one of the $q$ colours. No colour-preserving automorphism can act non-trivially within a component or interchange the components, so $D'(T)=q$.

If the optimal colouring of $(T',r)$ is unique up to a root-preserving automorphism, any colouring of $T$ with $q$ colours either fails to distinguish one of the rooted components or gives them isomorphic colourings. In the latter case, a colour-preserving automorphism interchanges the components. Hence $D'(T)>q$. One additional colour suffices: give both components optimal colourings, then recolour one edge in one component with a new colour. The two component colourings remain distinguishing and are now non-isomorphic, so $D'(T)=q+1$. Substituting these two possibilities for $q$ gives the claimed formulas.
\end{proof}

The second case occurs for a broad family of rooted trees. Fix a rooted tree $(T_0,r_0)$ and an integer $k\geq D'((T_0,r_0))$, and let $a$ be the number of non-isomorphic distinguishing edge colourings of $(T_0,r_0)$ from the palette $[k]$. Attach the roots of $ka$ copies of $(T_0,r_0)$ to a new root $r$. Each copy can be distinguished from the others by the colour of its edge to $r$ and its internal colouring. There are exactly $ka$ such combinations, so a distinguishing colouring from $[k]$ must use each combination once. The resulting colouring is unique up to a root-preserving automorphism, and fewer than $k$ colours cannot suffice. This construction leads to the following question, which concerns edge colourings rather than orientations.

\textbf{Question.} Characterise the rooted trees $(T',r)$ with a unique optimal distinguishing edge colouring up to a root-preserving automorphism.

\section{Graphs with $D'(G)=2$} \label{sec:D'=2}

If two edge colours suffice to break all non-trivial automorphisms, can edge directions alone do the same? In this section, we answer this question for traceable graphs and for connected claw-free graphs of order at least six.

A graph is \emph{traceable} if it has a Hamiltonian path. Pil\'sniak \cite{Pilsniak_traceable} proved that every traceable graph $G$ of order at least seven satisfies $D'(G)\leq 2$. The following construction gives a rigid orientation without any restriction on the order.

\begin{theorem} \label{thm: traceable}
If $G$ is traceable, then $OD'^-(G)=1$.
\end{theorem}
\begin{proof}
Let $v_1,\ldots,v_n$ be the vertices of a Hamiltonian path in order, and orient every edge $v_iv_j$ from $v_i$ to $v_j$ whenever $i<j$. The vertices reachable from $v_i$ by a non-empty directed path are exactly $v_{i+1},\ldots,v_n$. Thus $v_i$ reaches precisely $n-i$ vertices. These numbers are distinct and preserved by automorphisms, so every vertex is fixed.
\end{proof}

We next consider which automorphisms can preserve an orientation. An automorphism $\varphi\in\Aut(G)$ is \emph{twisted} if some positive power of $\varphi$ interchanges the endpoints of an edge, and \emph{non-twisted} otherwise. This distinction determines $OD'^+(G)$ when $D'(G)=2$.

\begin{theorem}
Let $G$ be a graph with $D'(G)=2$. If $G$ has a non-trivial, non-twisted automorphism, then $OD'^+(G)=2$. Otherwise, $OD'^-(G)=OD'^+(G)=1$.
\end{theorem}
\begin{proof}
A twisted automorphism cannot preserve any orientation of $G$. Indeed, if it did, all its powers would preserve that orientation, whereas one of them interchanges the endpoints of an edge and reverses its direction. Hence, if every non-trivial automorphism of $G$ is twisted, every orientation is rigid and $OD'^-(G)=OD'^+(G)=1$.

Now let $\varphi$ be a non-trivial, non-twisted automorphism. Consider the set
\[
A(G)=\{(u,v):uv\in E(G)\}
\]
of ordered pairs corresponding to the two possible directions of each edge. The map $\varphi$ induces a permutation $\varphi'$ of $A(G)$ given by $\varphi'(u,v)=(\varphi(u),\varphi(v))$. Since $\varphi$ is non-twisted, $(u,v)$ and $(v,u)$ lie in distinct cycles of $\varphi'$. Reversing every ordered pair in one cycle gives another cycle, so the cycles of $\varphi'$ occur in \emph{mirror pairs}.

Choose one cycle from each mirror pair and take all its ordered pairs as arcs. This selects exactly one direction for each edge, and $\varphi$ preserves the resulting orientation $\orG$. Since $\varphi$ is non-trivial, $D'(\orG)\geq 2$. On the other hand, every distinguishing edge colouring of $G$ also distinguishes $\orG$, so $D'(\orG)\leq D'(G)=2$. Thus $OD'^+(G)=2$.
\end{proof}

We now turn to claw-free graphs. Gorzkowska et al. \cite{GorzkowskaKargul} proved that every connected claw-free graph $G$ of order at least six satisfies $D'(G)\leq 2$. We adapt their greedy colouring procedure to construct a rigid orientation.

A \emph{path cover} of a graph $G$ is a collection of pairwise vertex-disjoint paths whose vertex sets partition $V(G)$. For each path, we choose one endpoint as its \emph{first vertex}.

We use the following two lemmas, stated as Lemma 5 and Claim 13 in \cite{GorzkowskaKargul}.

\begin{lemma}[\cite{GorzkowskaKargul}] \label{lem: pathcover}
Let $G$ be a connected claw-free graph and let $xy\in E(G)$. If $A\subset N(x)$ and $B\subset N(x)\setminus N[y]$, then:
\begin{enumerate}[1.]
\item $G[A]$ has a path cover with at most two paths.
\item $G[B]$ has a path cover with one path.
\end{enumerate}
\end{lemma}

\begin{lemma}[\cite{GorzkowskaKargul}] \label{lem: NsinC}
Let $G$ be a connected claw-free graph of order at least six, and let $C$ be a longest cycle of $G$. Then there is a vertex $s\in V(C)$ such that $N(s)\subseteq V(C)$.
\end{lemma}

In the 2-connected case, the construction starts by orienting the subgraph induced by a longest cycle. The following lemma handles this step when the cycle has a chord.

\begin{lemma} \label{lem: strong_rigid_hamiltonian}
Let $G$ be a Hamiltonian graph. Then $G$ admits a strongly connected rigid orientation if and only if $G$ is not a cycle. Moreover, if $G$ is not a cycle, then for every Hamiltonian cycle $C$ of $G$ there exists such an orientation in which $C$ is the unique directed Hamiltonian cycle.
\end{lemma}
\begin{proof}
If $G$ is a cycle, its only strongly connected orientations are directed cycles, which admit non-trivial rotations and are therefore not rigid.

Now suppose that $G$ is not a cycle, and fix a Hamiltonian cycle $C$ of $G$. Since $G$ has an edge outside $C$, choose a chord of $C$ and label the $n$ vertices in cyclic order as $w_0,w_1,\ldots,w_{n-1}$ so that the chosen chord is $w_0w_d$, where $2\leq d\leq n/2$. This is possible by choosing the starting endpoint and the direction of the labelling appropriately. Orient $C$ as $w_0\to w_1\to\cdots\to w_{n-1}\to w_0$. Orient the chosen chord from $w_d$ to $w_0$, and every other chord from its smaller-indexed endpoint to its larger-indexed endpoint.

To see that $C$ is the unique directed Hamiltonian cycle in this orientation, note that the only outgoing arc from $w_{n-1}$ is $w_{n-1}\to w_0$. Every directed Hamiltonian cycle must therefore contain this arc and cannot contain $w_d\to w_0$. Deleting $w_{n-1}\to w_0$ from such a cycle leaves a Hamiltonian path from $w_0$ to $w_{n-1}$ on which all indices strictly increase. The only possible vertex order is $w_0,w_1,\ldots,w_{n-1}$, so the cycle is $C$.

Consequently, every automorphism of this orientation restricts to a rotation of $C$. Suppose that a non-trivial rotation preserves the orientation. The subgroup it generates contains a rotation $w_i\mapsto w_{i+q}$, with indices taken modulo $n$, for some divisor $q$ of $n$ satisfying $1\leq q\leq n/2$. If $d+q<n$, this rotation maps the chord $w_d\to w_0$ to a different chord $w_{d+q}\to w_q$. This is impossible, since $w_d\to w_0$ is the only chord directed from a larger index to a smaller one. Otherwise, $d=q=n/2$, and the rotation reverses the chosen chord, which is also impossible. Thus the orientation of $G$ is rigid. It is also strongly connected, since it contains the directed Hamiltonian cycle $C$.
\end{proof}

We now apply these lemmas to construct rigid orientations of claw-free graphs.

\begin{theorem}
If $G$ is a connected claw-free graph of order at least six, then $OD'^-(G)=1$.
\end{theorem}
\begin{proof}
\emph{Case 1: $G$ is 2-connected.}
Let $C$ be a longest cycle of $G$. If $C$ is Hamiltonian, Theorem \ref{thm: traceable} gives the result. Otherwise, choose a vertex $u$ outside $C$ with a neighbour $v$ on $C$. The cycle has length at least four. By maximality of its length, $u$ is not adjacent to either neighbour of $v$ on $C$. Since $G$ is claw-free, these two neighbours must be adjacent to each other. Thus $C$ has a chord.

By Lemma \ref{lem: strong_rigid_hamiltonian}, orient $G[V(C)]$ so that it is rigid and $C$ is its unique directed Hamiltonian cycle. Choose $v_1\in V(C)$ with $N(v_1)\subseteq V(C)$, as guaranteed by Lemma \ref{lem: NsinC}, and label the remaining vertices $v_2,\ldots,v_n$ in the direction of $C$, where $n=|V(C)|$.

We extend the orientation using a set $R$ of reached vertices and a set $P\subseteq R$ of processed vertices. Initially, $R=V(C)$ and $P=\{v_1\}$. Every processed vertex has all its neighbours in $R$. We assign labels so that every vertex other than $v_1$ has a neighbour with a smaller label; this holds on $C$, and the rule below preserves it.

At each step, choose the vertex $v\in R\setminus P$ with the smallest label and put $S=N(v)\setminus R$. If $S$ is empty, add $v$ to $P$ and continue. Otherwise, $v$ has a neighbour $v'$ with a smaller label, which has already been processed. Since $N(v')\subseteq R$, we have
\[
S\subseteq N(v)\setminus N[v'].
\]
By Lemma \ref{lem: pathcover}, $G[S]$ has a Hamiltonian path. Label its vertices with the next consecutive integers in path order, and orient every edge of $G[S]$ from the smaller to the larger label. Orient all edges between $v$ and $S$ towards $S$, then replace $R$ by $R\cup S$ and add $v$ to $P$. The first vertex of the path has the smaller-labelled neighbour $v$, and every later vertex has its predecessor on the path, so the labelling condition is preserved.

Repeat until $R\setminus P$ is empty. At that point every vertex of $R$ has all its neighbours in $R$, so connectedness gives $R=P=V(G)$. Orient each remaining edge from the smaller to the larger label. All arcs between $V(C)$ and its complement now point away from $C$, and every arc with both endpoints outside $C$ increases the label. Hence every directed cycle lies in $V(C)$, and $C$ remains the unique directed cycle of length $n$.

We now prove that the completed orientation is rigid. Every automorphism preserves $C$ setwise and therefore restricts to an automorphism of the rigid orientation of $G[V(C)]$. Thus the initial set $R=V(C)$ is pointwise fixed. Suppose inductively that $R$ is pointwise fixed at the start of a step. Then $v$ is fixed and $S=N(v)\setminus R$ is setwise fixed. Its induced orientation is rigid by the argument in the proof of Theorem \ref{thm: traceable}: the vertices have distinct numbers of reachable vertices within $S$. Thus every vertex added to $R$ is fixed. Induction gives that all vertices of $G$ are fixed.

\emph{Case 2: $G$ is not 2-connected.}
Choose an end-block $B$ of $G$, with unique cut-vertex $v$; here a block may be a bridge $K_2$. Such a block exists because the block-cut tree has a leaf. Choose a neighbour $u$ of $v$ in $B$ and put $H=G-u$. Since $u$ is not a cut-vertex, $H$ is connected, and it is claw-free as an induced subgraph of $G$. We first orient $H$, leaving all edges incident to $u$ unoriented.

Partition $N_H(v)$ into
\[
S_1=N_H(v)\cap V(B),\qquad S_2=N_H(v)\setminus V(B).
\]
There are no edges between these sets, and $S_2$ is non-empty. Choose $w\in S_2$. Applying Lemma \ref{lem: pathcover} in $G$ to the edges $vw$ and $vu$ gives Hamiltonian paths in $G[S_1]$ and $G[S_2]$, respectively, whenever the corresponding set is non-empty. Label $v$ first, then the vertices of $S_1$ and $S_2$ in their respective path orders. Orient every edge within either set from the smaller to the larger label, and every edge from $v$ to $N_H(v)$ away from $v$.

Set $R=N_H[v]$ and $P=\{v\}$, and apply the recursive procedure from Case 1 entirely within $H$, taking all neighbourhoods in $H$. The initial labels satisfy the required condition, and every neighbour of the processed vertex $v$ is already in $R$. Complete the orientation of $H$ by directing the remaining edges from the smaller to the larger label. Finally, orient every edge $ux$ of $G$ from $u$ to $x$.

The vertex $u$ is the unique source. Indeed, $v$ receives an arc from $u$, and every other vertex of $H$ receives an incoming arc either from $v$ at the start or from the vertex that first reaches it. Hence $u$ is fixed, and so is $v$, the only cut-vertex of $G$ adjacent to $u$. Moreover, $B$ is preserved because it is the unique block containing $u$. It follows that $S_1$ and $S_2$ are each setwise fixed. Their induced orientations are rigid by the Hamiltonian path construction, so the initial set $R=N_H[v]$ is pointwise fixed. The same induction as in Case 1 now fixes every vertex of $H$, completing the proof.
\end{proof}

\section{Conclusions}

We have obtained exact formulas for $OD'^-(G)$ and $OD'^+(G)$ in terms of $D'(G)$ for bipartite graphs whose classes cannot be interchanged by an automorphism, and for trees. For arbitrary connected bipartite graphs of order at least three, Theorem \ref{thm: bipartite_bounds} restricts each parameter to at most two consecutive integer values. We have also constructed rigid orientations of traceable graphs and of connected claw-free graphs of order at least six. For Hamiltonian graphs other than cycles, our construction gives a strongly connected rigid orientation with a prescribed cycle as its unique directed Hamiltonian cycle.

The factor of two in our bounds reflects the use of arc directions to encode one of two edge colours. Our results suggest the following lower bound for arbitrary connected graphs.

\begin{conjecture}
If $G$ is a connected graph, then $OD'^-(G) \geq \lfloor D'(G)/2 \rfloor$.
\end{conjecture}

Theorem \ref{thm: bipartite_bounds} establishes this bound for all connected bipartite graphs of order at least three. The conjecture would extend this consequence beyond bipartite graphs. It would also imply that any  graph with a rigid orientation has distinguishing index at most three.

A complementary question is whether two colours in a distinguishing edge colouring can always be replaced by edge directions alone.

\begin{conjecture}
If $G$ is a connected graph with $D'(G)=2$, then $OD'^-(G) = 1$.
\end{conjecture}

Our results verify this conjecture for trees and traceable graphs, and for connected claw-free graphs of order at least six.

For balanced bipartite graphs whose classes can be interchanged by an automorphism, it remains to characterise when each of the two candidate values in Theorem \ref{thm: bipartite_bounds} is attained. Our result for trees gives such a characterisation in one class, but the general case remains open.


\begin{thebibliography}{99}
\bibitem{AlbertsonCollins}
M.~O.~Albertson, K.~L.~Collins, \emph{Symmetry breaking in graphs}, Electron. J. Combin., 3(1)(1996) \#18.
\bibitem{Babai}
L.~Babai, \emph{Asymmetric trees with two prescribed valences}, Acta Math. Acad. Sci. Hung. 29, (1977).
\bibitem{FisherIsaak}
M.~Fisher, G.~Isaak, \emph{Distinguishing colorings of Cartesian products of complete graphs}, Discrete Math. 308 (2008) 2240-2246.
\bibitem{GorzkowskaKargul}
A.~Gorzkowska, E.~Kargul, S.~Musia{\l}, K.~Pal, \emph{Edge-distinguishing of star-free graphs},  Electron. J. Combin. 27(3) (2020) \#P3.30.
\bibitem{ijk}
  W.~Imrich, J.~Jerebic and S.~Klav{\v{z}}ar,
  \emph{The Distinguishing Number of Cartesian Products of Complete Graphs},
  European J. Combin. 29 (2008) 922-929.
\bibitem{KalinowskiPilsniak}
R.~Kalinowski, M.~Pil\'sniak, \emph{Distinguishing graphs by edge-colourings}, European J. Combin. 45 (2015) 124-131.
\bibitem{KalinowskiPilsniak2}
R.~Kalinowski, M.~Pil\'sniak, \emph{Proper distinguishing arc-colourings of symmetric digraphs}, Appl. Math. Comput. 421 (2022)  art. no. 126939.
\bibitem{KalinowskiPilsniakProrok}
R.~Kalinowski, M.~Pil\'sniak, M.~Prorok, \emph{Distinguishing arc-colourings of symmetric digraphs}, Art Discrete Appl. Math. 6 (2023) \#P2.04.
\bibitem{KwasnyStawiski_regular}
J.~Kwa\'sny, M.~Stawiski, \emph{Distinguishing regular graphs}, preprint, arXiv:2207.14728.
\bibitem{MeslemSopena}
K.~Meslem, E.~Sopena, \emph{Distinguishing numbers and distinguishing indices of oriented graphs}, Discrete Appl. Math., 285 (2020) 330-342.
\bibitem{Pilsniak_traceable}
M.~Pil\'sniak, \emph{Improving upper bounds for the distinguishing index}, Ars Math. Contemp. 13 (2017) 259-274.
\bibitem{west}
D.B.~West, \emph{Introduction to Graph Theory}, 
Prentice Hall, Inc., 2nd Edition, 2001


\end{thebibliography}
\end{document}